\documentclass{elsart}
\usepackage{amssymb, amsfonts, amsmath, url, graphicx}

\def\a{{\alpha}}
\def\b{{\beta}}

\def\barr{\begin{array}}
\def\earr{\end{array}}

\def\0{\leqno}
\def\ov{\overline}

\def\Z{{\rlap{$\kern2pt{\rm Z}$}{\rm Z}\,}}

\def\fric{\displaystyle\frac{1}{2}}
\def\alege#1#2{{\displaystyle{{#1}\choose{#2}}}}
\def\n{\noindent }
\def\dd{\displaystyle}

\def\dd{\displaystyle}
\def\1{^{-1}}
\def\barr{\begin{array}}
\def\earr{\end{array}}

\def\a{{\alpha}}

\def\<{\left<}
\def\>{\right>}
\def\calf{{\cal F}}
\def\calc{{\cal C}}

\begin{document}
\begin{frontmatter}

\title{A new equivalence relation to classify\\ the fuzzy subgroups of finite groups}
\author{Marius T\u arn\u auceanu}
\ead{tarnauc@uaic.ro}
\address{Faculty of Mathematics,
\textquotedblleft Al.I. Cuza\textquotedblright\ University, Ia\c
si, Romania}

\maketitle

\begin{abstract}
In this paper a new equivalence relation $\approx$ to classify the
fuzzy subgroups of finite groups is introduced and studied. This
generalizes the equivalence relation $\sim$ defined on the lattice
of fuzzy subgroups of a finite group that has been used in our
previous papers (see e.g. \cite{16,24}). Explicit formulas for the
number of distinct fuzzy subgroups with respect to $\approx$ are
obtained in some particular cases.
\end{abstract}

\begin{keyword}
equivalence relations, fuzzy subgroups, chains of subgroups,
groups of automorphisms, group actions, fixed points, cyclic
groups, elementary abelian $p$-groups, dihedral groups, symmetric
groups.

{\it 2010 MSC:} Primary 20N25, 03E72; Secondary 05E15, 20D30,
20D60
\end{keyword}
\end{frontmatter}

\vspace{-10mm}
\section{Introduction}
\vspace{-5mm}
One of the most important problems of fuzzy group theory is to
classify the fuzzy subgroups of a finite group. This topic has
enjoyed a rapid development in the last few years. Several papers
have treated the particular case of finite abelian groups. Thus,
in \cite{8} the number of distinct fuzzy subgroups of a finite
cyclic group of square-free order is determined, while \cite{9},
\cite{10}, \cite{11} and \cite{27} deal with this number for
cyclic groups of order $p^n q^m$ ($p$, $q$ primes). Recall here
the paper \cite{24} (see also \cite{23}), where a recurrence
relation is indicated which can successfully be used to count the
number of distinct fuzzy subgroups for two classes of finite
abelian groups: cyclic groups and elementary abelian $p$-groups.
The explicit formula obtained for the first class leads in
\cite{16} to an important combinatorial result: a precise
expression of the well-known central Delannoy numbers in an
arbitrary dimension. Next, the study has been extended to some
remarkable classes of nonabelian groups: dihedral groups,
symmetric groups, finite $p$-groups having a cyclic maximal
subgroup and hamiltonian groups (see \cite{18}, \cite{20},
\cite{21} and \cite{25}, respectively). The same problems were
also investigated for the fuzzy normal subgroups of a finite group
(see \cite{22}).

Note that in all our papers mentioned above the fuzzy (normal)
subgroups of finite groups have been classified up to the same
natural equivalence relation $\sim$ defined on the fuzzy (normal)
subgroup lattices. This extends the equivalence relation used in
Murali's papers \cite{7}-\cite{11} and gives a powerful connection
between the fuzzy subgroups and certain chains of subgroups of
finite groups. Recall here the technique initiated in \cite{2}
(see also \cite{28}) to derive fuzzy theorems from their crisp
versions. Some other different approaches to classify the fuzzy
subgroups can be found in \cite{4} and \cite{5}.

In the present paper we will treat the problem of classifying the
fuzzy subgroups of a finite group $G$ by using a new equivalence
relation $\approx$ on the lattice $FL(G)$ of all fuzzy subgroups
of $G$. This is more general as $\sim$, excepting the case when
$G$ is cyclic (for which we will prove that
$\approx\hspace{1mm}=\hspace{1mm}\sim$). On the other hand, its
definition has a consistent group theoretical foundation, by
invol\-ving the knowledge of the automorphism group associated to
$G$. In order to count the distinct equivalence classes relative
to $\approx$, we shall use an interesting result of combinatorial
group theory: the Burnside's lemma (see \cite{13} or \cite{26}).
Our method will be exemplified for several remarkable classes of
finite groups. Also, we will compare the explicit formulas for the
numbers of distinct fuzzy subgroups with respect to $\approx$ with
the similar ones obtained in the case of $\sim$. Our approach is motivated
by the realization that in a theoretical study of fuzzy groups, fuzzy
subgroups are distinguished by their level subgroups and not by their images in
$[0,1]$. Consequently, the study of some equivalence relations between the chains 
of level subgroups of fuzzy groups is very important. It can also lead to other significant 
results which are similar with the analogous results in classical group theory. 

The paper is organized as follows. In Section 2 we present some
preliminary results on the fuzzy subgroups and the group actions
of a finite group $G$. Section 3 deals with a detailed description
of the new equivalence relation $\approx$ defined on $FL(G)$ and
of the technique that will be used to classify the fuzzy subgroups
of $G$. These are counted in Section 4 for the following classes
of finite groups: cyclic groups, elementary abelian $p$-groups,
dihedral groups and symmetric groups. In the final section several
conclusions and further research directions are indicated.

Most of our notation is standard and will usually not be repeated
here. Basic notions and results on lattices, groups and fuzzy
groups can be found in \cite{1}, \cite{14} and \cite{3} (see also
\cite{6}), respectively. For subgroup lattice concepts we refer
the reader to \cite{12} and \cite{15}.\vspace{-5mm}

\section{Preliminaries}
\vspace{-5mm}

Let $G$ be a group, $\calf(G)$ be the collection of all fuzzy subsets of $G$ and
$FL(G)$ be the lattice of fuzzy subgroups of $G$ (see e.g. \cite{6}). The fuzzy
(normal) subgroups of $G$ can be classified up to some
natural equiva\-lence relations on $\calf(G)$. One of them (used
in \cite{16}-\cite{25}, as well as in \cite{27}) is defined by
$$\mu\sim\eta\mbox{\ \ iff\ \ }(\mu(x)>\mu(y)\Longleftrightarrow \eta(x)>\eta(y),\mbox{ for all }x,y\in G)$$
and two fuzzy (normal) subgroups $\mu$, $\eta$ of  $G$ are said to
be {\it distinct} if $\mu\not\sim\eta$. This equivalence relation
generalizes that used in Murali's papers \cite{7}-\cite{11}. Also,
it can be connected to the concept of level subgroup. In this way,
suppose that the group $G$ is finite and let $\mu:G\to[0,1]$ be a
fuzzy (normal) subgroup of $G$. Put
$\mu(G){=}\{\alpha_1,\alpha_2,...,\alpha_r\}$ and assume that
\mbox{$\alpha_1>\alpha_2>...>\alpha_r.$} Then $\mu$ determines the
following chain of (normal) subgroups of $G$ which ends in $G$:
$$_\mu G_{\alpha_1}\subset {}_\mu G_{\alpha_2}\subset ...\subset {}_\mu G_{\alpha_r}=G.\vspace{-4mm}\leqno(*)$$

A necessary and sufficient condition for two fuzzy (normal)
subgroups $\mu,\eta$ of $G$ to be equivalent with respect to
$\sim$ has been identified in \cite{27}: $\mu\sim\eta$ if and only
if  $\mu$ and $\eta$ have the same set of level subgroups, that is
they determine the same chain of (normal) subgroups of type $(*)$.
This result shows that {\it there exists a bijection between the
equivalence classes of fuzzy {\rm (}normal\,{\rm )} subgroups of
$G$ and the set of chains of {\rm (}normal\,{\rm )} subgroups of
$G$ which end in $G$}. So, the problem of counting all distinct
fuzzy (normal) subgroups of $G$ can be translated into a
combinatorial problem on the subgroup lattice $L(G)$ (or on the
normal subgroup lattice $N(G)$) of $G$: finding the number of all
chains of (normal) subgroups of $G$ that terminate in $G$. Notice
also that in our previous papers we have denoted these numbers by
$h(G)$ (respectively by $h'(G)$).

Even for some particular classes of finite groups $G$, as finite
abelian groups, the problem of determining $h(G)$ is very
difficult. The largest classes of groups for which it was
completely solved are constituted by finite cyclic groups (see
Corollary 4 of \cite{24}) and by finite elementary abelian
$p$-groups (see the main result of \cite{23}). Recall that if
$\mathbb{Z}_n$ is the finite cyclic group of order $n$ and
$n=p_1^{m_1}p_2^{m_2}...p_s^{m_s}$ is the decomposition of $n$ as
a product of prime factors, then we have\vspace{-2mm}
$$h(\mathbb{Z}_n)=2^{\sum\limits_{\a=1}^s m_\a}\sum_{i_2=0}^{m_2}\ \sum_{i_3=0}^{m_3}...\sum_{i_s=0}^{m_s}\left(-\fric\right)^{\sum\limits_{\a=2}^s i_\a}\prod_{\a=2}^s\alege{m_\a}{i_\a}\alege{m_1+\sum\limits_{\b=2}^\a(m_\b-i_\b)}{m_\a},$$
where the above iterated sums are equal to 1 for $s=1$. The number
$h(G)$ has been also computed in several particular situations for
dihedral groups and symmetric groups in \cite{18} and \cite{20},
respectively. A complete study of finding the number $h'(G)$ for
all above classes of finite groups can be found in \cite{22}.

Finally, we recall that an {\it action} of the group $G$ on a
nonempty set $X$ is a map $\rho: X\times G\longrightarrow X$
satisfying the following two conditions:\vspace{-2mm}
\begin{enumerate}
\item[a)] $\rho(x,g_1g_2)=\rho(\rho(x,g_1),g_2)$, for all $g_1, g_2\in G$ and $x\in X$;
\item[b)] $\rho(x,e)=x$, for all $x\in X$.
\end{enumerate}\vspace{-2mm}
Every action $\rho$ of $G$ on $X$ induces an equivalence relation
$R_{\rho}$ on $X$, defined by
$$xR_{\rho}y \mbox{ if and only if there exists } g\in G \mbox{ such that }
y=\rho(x,g).\vspace{2mm}$$The equivalence classes of $X$ modulo $R_{\rho}$ are
called the {\it orbits} of $X$ relative to the action $\rho$. For
any $g \in G$, we denote by $Fix_{X}(g)$ the set of all elements
of $X$ which are fixed by $g$, that is
$$Fix_{X}(g)=\{x\in X \mid \rho(x,g)=x\}.$$If both $G$ and $X$ are finite,
then the number of distinct orbits of $X$ relative to $\rho$ is
given by the equality:\vspace{-2mm}
$$N=\frac{1}{|G|}\sum_{g\in G}|Fix_{X}(g)|.\vspace{-2mm}$$This result is known as the {\it
Burnside's lemma} and plays an important role in solving many
problems of combinatorics, finite group theory or graph theory. In
the next section it will be also used to compute the number of
distinct fuzzy subgroups of a finite group $G$ with respect to a
certain equivalence relation on $FL(G)$, induced by an action of
the automorphism group ${\rm Aut}(G)$ associated to $G$ on
$FL(G)$.\vspace{-5mm}

\section{A new equivalence relation on $FL(G)$}
\vspace{-5mm}
Let $G$ be a finite group. Then it is well-defined the following
action of ${\rm Aut}(G)$ on $FL(G)$
$$\hspace{-15mm}\rho: FL(G)\times {\rm Aut}(G)\longrightarrow FL(G),$$ $$\hspace{14mm}\rho(\mu, f)=\mu \circ f, \mbox{ for all } (\mu,f)\in FL(G)\times {\rm Aut}(G),\vspace{2mm}$$
and let us denote by $\approx_{\rho}$ the equivalence relation on
$FL(G)$ induced by $\rho$, namely
$$\mu\approx_{\rho}\eta \mbox{ if and only if there exists } f\in {\rm Aut}(G) \mbox{ such that }\eta = \mu \circ f.\vspace{-4mm}$$

As we have seen in Section 2, every fuzzy subgroup of $G$
determines a chain of subgroups of $G$ which ends in $G$ (that is,
a chain of type $(*)$). In this way, the above action $\rho$ can
be described in terms of chains of subgroups of $G$. Let $\mu,
\eta \in FL(G)$ and put
$\mu(G)=\{\alpha_1,\alpha_2,...,\alpha_m\}$ (where
$\alpha_1>\alpha_2>...>\alpha_m$),
$\eta(G)=\{\beta_1,\beta_2,...,\beta_n\}$ (where
$\beta_1>\beta_2>...>\beta_n$). If\vspace{-2mm}
$$\calc_{\mu}:\hspace{3mm}_\mu G_{\alpha_1}\subset {}_\mu
G_{\alpha_2}\subset ...\subset {}_\mu G_{\alpha_m}=G$$and
$$\hspace{-2mm}\calc_{\eta}:\hspace{3mm}_\eta G_{\beta_1}\subset
{}_\eta G_{\beta_2}\subset ...\subset {}_\eta G_{\beta_n}=G$$are
the corresponding chains of type $(*)$, then we can easily see
that $\mu\approx_{\rho}\eta$ if and only if $Im(\mu)=Im(\eta)$
(i.e. $m=n$ and $\alpha_i=\beta_i$, $i=\ov{1,n}$) and there is $f
\in {\rm Aut}(G)$ such that $f(_\eta G_{\beta_i})={_\mu
G_{\alpha_i}}$, for all $i=\ov{1,n}$. In other words, two fuzzy
subgroups $\mu$ and $\eta$ of $G$ are equivalent with respect to
$\approx_{\rho}$ if and only if they have the same image and there
is an automorphism $f$ of $G$ which maps $\calc_{\eta}$ into
$\calc_{\mu}$.

Since we are interested to work only with chains of subgroups, in
the following we will consider the equivalence relation $\approx$
on $FL(G)$ defined by
$$\mu\approx\eta \mbox{ iff } \exists f\in {\rm Aut}(G) \mbox{ such that } f(\calc_{\eta})=\calc_{\mu}.$$
Obviously, this is a little more general as $\approx_{\rho}$: we
can easily prove that if $\mu\approx\eta$, then their images are
not necessarily equal, but certainly there is a bijection between
$Im(\mu)$ and $Im(\eta)$. Moreover, we also remark that $\approx$
generalizes the equivalence relation $\sim$ (recall that
$\mu\sim\eta \Longleftrightarrow \calc_{\mu}=\calc_{\eta}$, i.e.
the above automorphism $f$ is the identical map of $G$).

Next, we will focus on computing the number $N$ of distinct fuzzy
subgroups of $G$ with respect to $\approx$, that is the number of
distinct equivalence classes of $FL(G)$ modulo $\approx$. Denote
by $\ov{\calc}$ the set consisting of all chains of subgroups of
$G$ terminated in $G$. Then the previous action $\rho$ of ${\rm
Aut}(G)$ on $FL(G)$ can be seen as an action of ${\rm Aut}(G)$ on
$\ov{\calc}$ and $\approx_{\rho}$ as the equivalence relation
induced by this action. An equivalence relation on $\ov{\calc}$
which is similar with $\approx$ can also be constructed in the
following manner: for two chains
$$\calc_1:\hspace{3mm} H_1 \subset H_2 \subset ...\subset H_m=G\,\,\, \mbox{ and }\,\,\,
\calc_2:\hspace{3mm} K_1 \subset K_2 \subset ...\subset
K_n=G$$of $\ov{\calc}$, we put $$\calc_1 \approx \calc_2 \mbox{
iff } m=n \mbox{ and }\exists f \in {\rm Aut}(G) \mbox{ such that
} f(H_i)=K_i, i=\ov{1,n}.$$In this case the orbit of a chain
$\calc \in \ov{\calc}$ is $\{f(\calc) \mid f \in {\rm Aut}(G)\}$,
while the set of all chains in $\ov{\calc}$ that are fixed by an
automorphism $f$ of $G$ is $Fix_{\ov{\calc}}(f)=\{\calc \in \ov{\calc} \mid\newline f(\calc)=\calc \}$.
Now, the number $N$ is obtained by applying the Burnside's lemma:
$$N=\frac{1}{|{\rm Aut}(G)|}\sum_{f\in {\rm Aut}(G)}|Fix_{\ov{\calc}}\hspace{1mm}(f)|.\vspace{-4mm}\leqno(**)$$

Finally, we note that the above formula can successfully be used
to calculate $N$ for any finite group $G$ whose subgroup lattice
$L(G)$ and automorphism group ${\rm Aut}(G)$ are known.\vspace{-5mm}

\section{The number of distinct fuzzy subgroups of finite groups}
\vspace{-5mm}
In this section we shall compute explicitly the number $N$ of all
distinct fuzzy subgroups with respect to $\approx$ for several
remarkable classes of finite groups. The comparison with the
similar results relative to $\sim$ obtained in our previous papers
is also made.

\smallskip
{\bf 4.1. The number of distinct fuzzy subgroups of finite
cyclic groups}
\smallskip

Let $\mathbb{Z}_n=\{\,\ov{0}, \ov{1}, ..., \ov{n-1}\,\}$ be the
finite cyclic group of order $n$. The subgroup structure of
$\mathbb{Z}_n$ is well-known (see \cite{14}, I): for every divisor
$d$ of $n$, there is a unique subgroup of order $d$ in
$\mathbb{Z}_n$, namely $\langle\ov{\frac{n}{d}}\rangle$. Moreover,
the following lattice isomorphism holds $$L(\mathbb{Z}_n)\cong
L_n,$$where $L_n$ denotes the lattice of all divisors of $n$,
under the divisibility. It is also well-known the structure of the
automorphism group of $\mathbb{Z}_n$: if $f_{\ov{d}}: \mathbb{Z}_n
\longrightarrow \mathbb{Z}_n$ is the map defined by
$f_{\ov{d}}\hspace{1mm}(\ov{k})=\ov{dk}$, for all
$\ov{d},\ov{k}\in\mathbb{Z}_n$, then
$${\rm Aut}(\mathbb{Z}_n)=\{f_{\ov{d}} \hspace{1mm}\mid\hspace{1mm}
(d,n)=1\}.$$In particular, we have $|{\rm
Aut}(\mathbb{Z}_n)|=\varphi(n)$, where $\varphi$ is the Euler's
totient function.

In order to find the number $N$ for $\mathbb{Z}_n$, the equality
$(**)$ can be used. Ne\-ver\-theless, in this case we shall prefer to
give a direct solution, founded on the remark that the equivalence
relations $\approx$ and $\sim$ coincide for a finite cyclic group.

\n{\bf Theorem 4.1.1.} {\it For a finite group $G$, the following
conditions are equivalent:
\begin{enumerate}
\item[a{\rm )}] The equivalence relations $\approx$ and $\sim$ associated to $G$ coincide.
\item[b{\rm )}] $G$ is cyclic.
\end{enumerate}}

\n{\bf Proof.} Assume first that $G=\mathbb{Z}_n$ is cyclic. Let
$\calc_1, \calc_2$ be two chains of subgroups of
$\mathbb{Z}_n$ ended in $\mathbb{Z}_n$ ($\calc_1:\hspace{1mm} H_1
\subset H_2 \subset ...\subset H_p=\mathbb{Z}_n, \hspace{1mm}
\calc_2:\hspace{1mm} K_1 \subset K_2 \subset ...\subset
K_q=\mathbb{Z}_n$) which satisfy $\calc_1 \approx \calc_2$ and
take $f\in {\rm Aut}(\mathbb{Z}_n)$ such that $f(\calc_1)=
\calc_2$. Then $p=q$ and $f(H_i)=K_i$, for any $i=\ov{1,p}$. Since
$f$ is an automorphism, the subgroups $H_i$ and $K_i$ are of the
same order, therefore they must be equal. This shows that
$\calc_1=\calc_2$, that is for a finite cyclic group the
equivalence relations $\approx$ and $\sim$ are the same.

Conversely, let us assume that for a finite group $G$ we have
$\approx\hspace{1mm}=\hspace{1mm}\sim$ and take an arbitrary subgroup
$H$ of $G$. Then, for any $x\in G$, the chains $H \subset G$ and $H^x \subset G$
are equivalent modulo $\approx$, because there exists the inner automorphism $f_x$ of
$G$ ($f_x(g)=xgx^{-1}$, for all $g\in G$) such that $f_x(H)=H^x$.
This shows that $H=H^x$ for all $x\in G$, that is $H$ is a normal
subgroup of $G$. Thus $G$ is either a hamiltonian group or an
abelian group. In the first case $G$ is of type
$Q\times\mathbb{Z}_2^k \times A$, where $Q=\langle x, y \mid
x^4=1, x^2=y^2, y^{-1}xy=x^{-1}\rangle$ is the quaternion group of
order 8, $\mathbb{Z}_2^k$ is the direct product of $k$ copies of
$\mathbb{Z}_2$ and $A$ is a finite abelian group of odd order.
Clearly, the function $f:G\longrightarrow G$ that maps $x$ into
$y$ and leaves invariant every element of $\mathbb{Z}_2^k \times
A$ is an automorphism of $G$. Then the chains $\langle x\rangle
\subset G$ and $\langle y\rangle \subset G$ are equivalent modulo
$\approx$, but they are not equal, a contradiction. In the second
case $G$ is a direct product of type $G_1 \times G_2 \times
\cdot\cdot\cdot \times G_s$, where each $G_i$ is an abelian
$p_i$-group, $i=\ov{1,s}$. If we assume that there exists $i\in
\{1, 2, ..., s\}$ such that the group $G_i$ is not cyclic, then
$G_i$ possesses two distinct maximal subgroups $M_i^1, M_i^2$ and
an automorphism $f_i$ satisfying $f_i(M_i^1)=M_i^2$. It is easy
to see that $f_i$ can be extended to an automorphism $f$ of $G$
mapping $M_i^1$ into $M_i^2$. Again, by our hypothesis, the chains
$M_i^1 \subset G$ and $M_i^2 \subset G$ must be equal. This leads
to $M_i^1=M_i^2$, a contradiction. Therefore every $G_i$,
$i=\ov{1,s}$, is cyclic and so is $G$ itself.
\hfill\rule{1,5mm}{1,5mm}

By Theorem 4.1.1 one obtains that $N=h(\mathbb{Z}_n)$ and thus the explicit
formula of $h(\mathbb{Z}_n)$ presented in Section 2 can be used to compute $N$, too.
Hence the following result holds.

\n{\bf Theorem 4.1.2.} {\it Let $n=p_1^{m_1}p_2^{m_2}...p_s^{m_s}$ be the decomposition
of $n\in\mathbb{N}, n\geq 2,$ as a product of prime power factors. Then the
number $N$ of all distinct fuzzy subgroups with respect to
$\approx$ of the finite cyclic group $\mathbb{Z}_n$ is given by
the equality\vspace{-2mm}
$$N=2^{\sum\limits_{\a=1}^s m_\a}\sum_{i_2=0}^{m_2}\ \sum_{i_3=0}^{m_3}...\sum_{i_s=0}^{m_s}\left(-\fric\right)^{\sum\limits_{\a=2}^s i_\a}\prod_{\a=2}^s\alege{m_\a}{i_\a}\alege{m_1+\sum\limits_{\b=2}^\a(m_\b-i_\b)}{m_\a},$$
where the above iterated sums are equal to {\rm 1} for $s=1$.}

\smallskip
{\bf 4.2. The number of distinct fuzzy subgroups of finite
elementary abelian $p$-groups}
\smallskip

It is well-known (for example, see \cite{14}) that a finite
abelian group can be written as a direct product of $p$-groups.
Therefore the problem of counting the fuzzy subgroups of finite
abelian groups must be first solved for $p$-groups. A special
class of abelian $p$-groups is constituted by elementary abelian
$p$-groups. Such a group $G$ has a direct decomposition of type\vspace{-3mm}
$$\mathbb{Z}^n_p=\underbrace{\mathbb{Z}_p\times\mathbb{Z}_p\times\cdots\times\mathbb{Z}_p}_{n\ {\rm factors}}\,,\vspace{-3mm}$$
where $p$ is a prime and $n\in\mathbb{N}^*$, and its distinct
fuzzy subgroups with respect to $\sim$ have been counted in
\cite{23,24}. It is also known that $\mathbb{Z}^n_p$ possesses a
natural linear space structure over the field $\mathbb{Z}_p$ and
that the automorphisms of the group $\mathbb{Z}^n_p$ coincide with
the automorphisms of this linear space. In this way, the
automorphism group ${\rm Aut}(\mathbb{Z}^n_p)$ is isomorphic to
the general linear group $GL(n,p)$ and, in particular, we have\vspace{-3mm}
$$|{\rm Aut}(\mathbb{Z}^n_p)|=\prod_{i=0}^{n-1}\hspace{1mm}(p^n-p^i).\vspace{-3mm}$$
Let $f \in {\rm Aut}(\mathbb{Z}^n_p)$ and $\calc \in
Fix_{\ov{\calc}}\hspace{1mm}(f)$, where $\calc: \hspace{1mm} H_1
\subset H_2 \subset ...\subset H_m=\mathbb{Z}^n_p$. Then
$f(\calc)=\calc$, that is $f(H_i)=H_i$, for all $i=\ov{1,m}$. This
shows that all subspaces $H_i$, $i=1,2, ..., m$, of
$\mathbb{Z}^n_p$ are invariant with respect to $f$. So, the
problem of counting the number of elements of
$Fix_{\ov{\calc}}\hspace{1mm}(f)$ is reduced to a linear algebra
problem: determining all chains of subspaces of the linear space
$\mathbb{Z}^n_p$ which ends in $\mathbb{Z}^n_p$ and are invariant
with respect to $f$.

Next, we will solve the above problem in the simplest case when
$p=n=2$ (note that the general case can be treated in a similar
manner). The proper subspaces of $\mathbb{Z}^2_2$, that in fact
correspond to the proper subgroups of the Klein's group, are
$H_1=\langle(\hat{1}, \hat{0})\rangle$, $H_2=\langle(\hat{0},
\hat{1})\rangle$ and $H_3=\langle(\hat{1}, \hat{1})\rangle$. As we
already have seen, every group automorphism of $\mathbb{Z}^2_2$ is
perfectly determined by a matrix contained in\vspace{-2mm}
$$\hspace{-15mm}GL(2,2)=\{A_1=\left(\begin{array}{cc}
\hat{1} & \hat{0} \\
\hat{0} & \hat{1}
\end{array}\right),\hspace{1mm} A_2=\left(\begin{array}{cc}
\hat{1} & \hat{0} \\
\hat{1} & \hat{1}
\end{array}\right),\hspace{1mm} A_3=\left(\begin{array}{cc}
\hat{0} & \hat{1} \\
\hat{1} & \hat{0}
\end{array}\right),\vspace{-2mm}$$ $$\hspace{19mm}A_4=\left(\begin{array}{cc}
\hat{1} & \hat{1} \\
\hat{0} & \hat{1}
\end{array}\right),\hspace{1mm} A_5=\left(\begin{array}{cc}
\hat{0} & \hat{1} \\
\hat{1} & \hat{1}
\end{array}\right),\hspace{1mm} A_6=\left(\begin{array}{cc}
\hat{1} & \hat{1} \\
\hat{1} & \hat{0}
\end{array}\right)\}\cong S_3.\vspace{-2mm}$$

For each $i\in \{1,2, ...,6\}$, let us denote by $f_i$ the
automorphism induced by the matrix $A_i$ and by
$L_{f_i}(\mathbb{Z}^2_2)$ the set consisting of all subspaces of
$\mathbb{Z}^2_2$ which are invariant relative to $f_i$, i.e.\vspace{-2mm}
$$L_{f_i}(\mathbb{Z}^2_2)=\{H
\leq_{\hspace{1mm}\mathbb{Z}_2}\hspace{-1mm}\mathbb{Z}^2_2 \mid
f_i(H)=H\}.\vspace{-2mm}$$We easily obtain $L_{f_1}(\mathbb{Z}^2_2)=\{1, H_1, H_2, H_3, \mathbb{Z}^2_2\}$,
$L_{f_2}(\mathbb{Z}^2_2)=\{1, H_2, \mathbb{Z}^2_2\}$, $L_{f_3}(\mathbb{Z}^2_2)\newline =\{1, H_3, \mathbb{Z}^2_2\}$,
$L_{f_4}(\mathbb{Z}^2_2)=\{1, H_1, \mathbb{Z}^2_2\}$,$L_{f_5}(\mathbb{Z}^2_2)=L_{f_6}(\mathbb{Z}^2_2)=\{1, \mathbb{Z}^2_2\}$,
where 1 denotes the trivial subspace of $\mathbb{Z}^2_2$. The above equalities show that $|Fix_{\ov{\calc}}\hspace{1mm}(f_1)|=8$,
$|Fix_{\ov{\calc}}\hspace{1mm}(f_i)|=4$ for $i=2,3,4$, and
$|Fix_{\ov{\calc}}\hspace{1mm}(f_i)|=2$ for $i=5,6$. Clearly, by
applying the formula $(**)$, we are now able to compute explicitly
the number $N$ associated to the elementary abelian 2-group
$\mathbb{Z}^2_2$.

\n{\bf Theorem 4.2.1.} {\it The number $N$ of all distinct fuzzy
subgroups with respect to $\approx$ of the elementary abelian {\rm
2}-group $\mathbb{Z}^2_2$ is given by the equality\vspace{-2mm}
$$N=\dd\frac{1}{6}\hspace{1mm}(8+3\cdot4+2\cdot2)=4.\vspace{-2mm}$$}
\n{\bf Remark.} By \cite{23,24} we know that
$h(\mathbb{Z}^2_2)=8$, i.e. there are 8 distinct chains of
subgroups (subspaces) of $\mathbb{Z}^2_2$ ended in
$\mathbb{Z}^2_2$, namely $\mathbb{Z}^2_2$,
$1\subset\mathbb{Z}^2_2$, $H_i\subset\mathbb{Z}^2_2$, $i=1,2,3$,
and $1\subset H_i\subset\mathbb{Z}^2_2$, $i=1,2,3$. From these
chains several are equivalent modulo $\approx$. More precisely,
the distinct equivalence classes with respect to $\approx$ of the
set $\ov{\calc}$ described above are: $\ov{\calc}_1=\{\mathbb{Z}^2_2\}$,
$\ov{\calc}_2=\{1\subset\mathbb{Z}^2_2\}$, $\ov{\calc}_3=\{H_1\subset\mathbb{Z}^2_2, H_2\subset\mathbb{Z}^2_2, H_3\subset\mathbb{Z}^2_2\}$,
$\ov{\calc}_4=\{1\subset H_1\subset\mathbb{Z}^2_2, 1\subset H_2\subset\mathbb{Z}^2_2, 1\subset H_3\subset\mathbb{Z}^2_2\}$.

\smallskip
{\bf 4.3. The number of distinct fuzzy subgroups of finite
dihedral groups}
\smallskip\vspace{-4mm}

The finite dihedral group $D_{2n}$ $(n\ge2)$ is the symmetry group
of a regular polygon with $n$ sides and has the order $2n$. The
most convenient abstract description of $D_{2n}$ is obtained by
using its generators: a rotation $a$ of order $n$ and a reflection
$b$ of order 2. Under these notations, we have
$$D_{2n}=\langle a, b \mid a^n=b^2=1, bab=a^{-1}\rangle=\{1, a, a^2, ..., a^{n-1}, b, ab, a^2b, ..., a^{n-1}b\}.$$The
automorphism group of $D_{2n}$ is well-known, namely
$${\rm Aut}(D_{2n})=\{f_{\alpha, \beta} \mid \alpha=\ov{0,n-1} \mbox{ with } (\alpha, n)=1,
\beta=\ov{0,n-1})\},$$where $f_{\alpha, \beta}: D_{2n} \longrightarrow D_{2n}$ is defined by
$f_{\alpha, \beta}\hspace{1mm}(a^i)=a^{\alpha i}$ and $f_{\alpha, \beta}\hspace{1mm}(a^ib)=a^{\alpha i+\beta}b$,
for all $i=\ov{0,n-1}$. We infer that
$$|{\rm Aut}(D_{2n})|=n\varphi(n).$$The structure of the subgroup
lattice $L(D_{2n})$ of $D_{2n}$ is also well-known: for every
divisor $r$ or $n$, $D_{2n}$ possesses a subgroup isomorphic to
$\mathbb{Z}_r$, namely $H^r_0=\hspace{1mm}\langle a^{\frac
nr}\rangle$, and $\frac{n}{r}$ subgroups isomorphic to $D_r$,
namely $H^r_i=\langle a^{\frac nr},a^{i-1}b\rangle,$
$i=1,2,...,\frac{n}{r}\hspace{1mm}.$

Next, for each $f_{\alpha, \beta} \in {\rm Aut}(D_{2n})$, let
$Fix(f_{\alpha, \beta})$ be the set consisting of all subgroups of
$D_{2n}$ that are invariant relative to $f_{\alpha, \beta}$, that
is
$$Fix(f_{\alpha, \beta})=\{H \leq D_{2n} \mid f_{\alpha,
\beta}\hspace{1mm}(H)=H\}.$$By using some elementary results of
group theory, these subsets of $L(D_{2n})$ can precisely be
determined: a subgroup of type $H^r_0$ belongs to $Fix(f_{\alpha,
\beta})$ if and only if $(\alpha, r)=1$, while a subgroup of type
$H^r_i$ belongs to $Fix(f_{\alpha, \beta})$ if and only if
$(\alpha, r)=1$ and $\frac{n}{r}$ divides $(\alpha-1)(i-1)+\beta$.

Under the above notation, computing
$|Fix_{\ov{\calc}}\hspace{1mm}(f_{\alpha, \beta})|$ is reduced to
computing the number of chains of $L(D_{2n})$ which end in
$D_{2n}$ and are contained in the set $Fix(f_{\alpha, \beta})$.
After then, an explicit expression of the number $N$ associated to
the group $D_{2n}$ will follow from $(**)$.

In the following we will apply this algorithm for the dihedral
group $D_8$. As we already have seen, the group ${\rm Aut}(D_8)$
consists of $4\varphi(4)=4\cdot2=8$ ele\-ments, namely ${\rm Aut}(D_8)=\{f_{1,0}, f_{1,1}, f_{1,2}, f_{1,3}, f_{3,0}, f_{3,1}, f_{3,2},
f_{3,3}\}$. The corresponding sets of subgroups of $D_8$ that are
invariant relative to the above automorphisms can be described by
a direct calculation: $Fix(f_{1,0}){=}L(D_8)$, $Fix(f_{1,1}){=}Fix(f_{1,3}){=}Fix(f_{3,1}){=}Fix(f_{3,3}){=}\{H^1_0, H^2_0, H^4_0, H^4_1\}$,
$Fix(f_{1,2}){=}\{H^1_0,\newline H^2_0,H^4_0, H^2_1, H^2_2, H^4_1\}$, $Fix(f_{3,0}){=}\{H^1_0, H^2_0, H^4_0, H^1_1, H^1_3, H^2_1, H^2_2, H^4_1\}$,
$Fix(f_{3,2}){=}\newline\{H^1_0, H^2_0, H^4_0, H^1_2, H^1_4, H^2_1, H^2_2, H^4_1\}$. Then we easily obtain:
$|Fix_{\ov{\calc}}\hspace{1mm}(f_{1,0})|=32$,\newline $|Fix_{\ov{\calc}}\hspace{1mm}(f_{1,1})|=|Fix_{\ov{\calc}}\hspace{1mm}(f_{1,3})|=|Fix_{\ov{\calc}}\hspace{1mm}(f_{3,1})|=|Fix_{\ov{\calc}}\hspace{1mm}(f_{3,3})|=8$,
$|Fix_{\ov{\calc}}\hspace{1mm}(f_{1,2})|=16$, $|Fix_{\ov{\calc}}\hspace{1mm}(f_{3,0})|=Fix_{\ov{\calc}}\hspace{1mm}(f_{3,2})|=24$. Thus the following result holds.

\n{\bf Theorem 4.3.1.} {\it The number $N$ of all distinct fuzzy
subgroups with respect to $\approx$ of the dihedral group $D_8$ is
given by the equality\vspace{-2mm}
$$N=\dd\frac{1}{8}\hspace{1mm}(32+4\cdot8+16+2\cdot24)=16.\vspace{-2mm}$$}
\n{\bf Remark.} In Theorem 3 of \cite{18} we have proved that the
number $h(D_{2p^m})$ of all distinct fuzzy subgroups relative to
$\sim$ of the dihedral group $D_{2p^m}$ (where $p$ is prime and
$m\in\mathbb{N^*}$) is\vspace{-2mm}
$$h(D_{2p^m})=\dd\frac{2^m}{p-1}(p^{m+1}+p-2).\vspace{-2mm}$$
In particular, we find $h(D_{2^m})=2^{2m-1}$ and so $h(D_8)=32$, a
number which is different from the number $N$ given by Theorem
4.3.1. We also remark that since $h(D_8)$ counts all chains of
subgroups of $D_8$ ended in $D_8$, it is in fact equal to the
number of these chains which are invariant relative to
$f_{1,0}=1_{D_8}$ (the identical automorphism of $D_8$), that is
to $|Fix_{\ov{\calc}}\hspace{1mm}(f_{1,0})|$.

\smallskip
{\bf 4.4. The number of distinct fuzzy subgroups of finite
symmetric groups}
\smallskip\vspace{-2mm}

In order to apply the formula $(**)$ for the symmetric group
$S_n$, $n\geq3$, we need to known the automorphism group ${\rm
Aut}(S_n)$. An important normal subgroup of this group is the
inner automorphism group ${\rm Inn}(S_n)$, which consists of all
automorphisms of $S_n$ of type $f_\sigma$, $\sigma \in S_n$, where
$f_\sigma(\tau)=\sigma\tau\sigma^{-1}$, for all $\tau \in S_n$. On
the other hand, we know that if $n\geq3$, then ${\rm Inn}(S_n)$ is
isomorphic to $S_n$, because the center $Z(S_n)$ of $S_n$ is
trivial and $S_n/Z(S_n)\cong {\rm Inn}(S_n)$. According to
\cite{14}, I, for any $n\neq6$ we have \vspace{-3mm}
$${\rm Aut}(S_n)={\rm Inn}(S_n)\cong S_n,\vspace{-3mm}$$while for $n=6$ we have\vspace{-2mm}
$${\rm Aut}(S_6)\cong {\rm Aut}(A_6) \mbox{ and }
({\rm Aut}(S_6):{\rm Inn}(S_6))=2.$$In particular, one obtains
$$|{\rm Aut}(S_n)|=\left\{\barr{lll}
n!\hspace{0,5mm},&n\neq6\\
&&\\
2\cdot6!\hspace{0,5mm},&n=6.\earr\right.$$ In the following we
will focus only on the case $n\neq6$. Then every automorphism of
$S_n$ is of the form $f_\sigma$ with $\sigma\in S_n$. Let $\calc
\in Fix_{\ov{\calc}}\hspace{1mm}(f_{\sigma})$, where $\calc:
\hspace{1mm} H_1 \subset H_2 \subset ...\subset H_m=S_n$. Then
$f_{\sigma}(\calc)=\calc$, that is $f_{\sigma}(H_i)=H_i$, for all
$i=\ov{1,m}$. This shows that $H_i^{\sigma}=H_i$, i.e. $\sigma$ is
contained in the normalizer $N_{S_n}(H_i)$ of $H_i$ in $S_n$,
$i=\ov{1,m}$. Therefore we have\vspace{-2mm}
$$\calc\in Fix_{\ov{\calc}}\hspace{1mm}(f_{\sigma})\Longleftrightarrow \sigma \in N_{S_n}(H_i), \mbox{ for all }
i=\ov{1,m} \Longleftrightarrow \sigma \in \dd\bigcap_{i=1}^m
\hspace{1mm}N_{S_n}(H_i),\0(1)$$which allows us to compute
explicitly $|Fix_{\ov{\calc}}\hspace{1mm}(f_{\sigma})|$. We will
exemplify this method for the symmetric group $S_3$. It is
well-known that $S_3$ has 6 elements, more precisely $$S_3=\{e,
\tau_1=(23), \tau_2=(13), \tau_3=(12), \sigma=(123),
\sigma^2=(132)\},$$and its subgroup lattice consists of the
trivial subgroup $H_0=\{e\}$, three subgroups of order 2, namely
$H_i=\langle\tau_i\rangle=\{e, \tau_i\}, \hspace{1mm}i=1,2,3$, a
subgroup of order 3, namely $H_4=\langle\sigma\rangle=\{e, \sigma,
\sigma^2\}$, and a subgroup of order 6, namely $H_5=S_3$. Then the
chains of subgroups of $S_3$ ended in $S_3$ are $H_5$, $H_i
\subset H_5$, $i=\ov{0,4}$, and $H_0 \subset H_i \subset H_5$,
$i=\ov{1,4}$. Also, we can easily find the normalizers of all
above subgroups: $N_{S_3}(H_0)=N_{S_3}(H_0)=N_{S_3}(H_0)=S_3$ (in
other words $H_0$, $H_4$ and $H_5$ are normal in $S_3$), and
$N_{S_3}(H_i)=H_i$, $i=1,2,3$. By using (1), one obtains
$|Fix_{\ov{\calc}}\hspace{1mm}(f_e)|=10$,
$|Fix_{\ov{\calc}}\hspace{1mm}(f_{\tau_i})|=6$, $i=1,2,3$, and
$|Fix_{\ov{\calc}}\hspace{1mm}(f_{\sigma})|=|Fix_{\ov{\calc}}\hspace{1mm}(f_{\sigma^2})|=4$.
So, we are able to compute explicitly the number $N$ associated to
the symmetric group $S_3$ (and also to the dihedral group $D_6$,
since $S_3\cong D_6$), in view of the formula $(**)$.

\n{\bf Theorem 4.4.1.} {\it The number $N$ of all distinct fuzzy
subgroups with respect to $\approx$ of the symmetric group $S_3$
is given by the equality\vspace{-2mm}
$$N=\dd\frac{1}{6}\hspace{1mm}(10+3\cdot6+2\cdot4)=6.\vspace{-2mm}$$}We remark that $N$ is different from $h(S_3)=10$
(computed in \cite{20}), as we expected. This is due to the fact
that, from the chains of $\ov{\calc}$ described above, $H_i
\subset H_5$, $i=1,2,3$, are equivalent modulo $\approx$, and the
same thing can be said about $H_0 \subset H_i \subset H_5$,
$i=1,2,3$. Finally, we remark that the above reasoning can
successfully be applied to compute the number $N$ associated to
the symmetric group $S_4$, whose subgroup structure has been
completely described in \cite{20}, too.\vspace{-8mm}

\section{Conclusions and further research}
\vspace{-5mm}
The study concerning the classification of the fuzzy subgroups of
(finite) groups is a significant aspect of fuzzy group theory. It
can be made with respect to some natural equivalence relations on
the fuzzy subgroup lattices, as the Murali's equivalence relation
used in \cite{7}-\cite{11}, $\sim$ used in \cite{16}-\cite{25} and
\cite{27}, or $\approx$ used in this paper. Obviously, other such
relations can be introduced and investigated. The problem of
counting the distinct fuzzy subgroups relative to the above
equivalence relations can also be extended to other remarkable
classes of finite groups. This will surely constitute the subject
of some further research.

Two open problems with respect to this topic are the
following.

\n{\bf Problem 1.} Generalize the results of Section 4, by
establishing some explicit formulas for the number of distinct
fuzzy subgroups of {\it arbitrary} elementary abelian $p$-groups,
dihedral groups and symmetric groups.

\n{\bf Problem 2.} Classify the fuzzy normal subgroups of a finite
group with respect to the equivalence relation $\approx$ defined
in Section 3. Use the particular classes of finite groups studied
in \cite{22}, whose normal subgroup structure can completely be
described.\vspace{5mm}

{\bf Acknowledgements.} The author is grateful to the reviewers for
their remarks which improve the previous version of the paper.


\begin{thebibliography}{00}
\bibitem{1} G. Gr\"atzer, General lattice theory, Academic Press, New York, 1978.
\bibitem{2} T. Head, A metatheorem for deriving fuzzy theorems from crisp versions, Fuzzy Sets and Systems {\bf 73} (1995), 349-358; {\bf 79} (1996), 277-278.
\bibitem{3} R. Kumar, Fuzzy algebra, I, Univ. of Delhi, Publ. Division, 1993.
\bibitem{4} M. Mashinchi, M. Mukaidono, A classification of fuzzy subgroups, Ninth Fuzzy System Symposium, Sapporo, Japan, 1992, 649-652.
\bibitem{5} M. Mashinchi, M. Mukaidono, On fuzzy subgroups classification, Research Report of Meiji Univ. {\bf 9} (1993), 31-36.
\bibitem{6} J.N. Mordeson, K.R. Bhutani, A. Rosenfeld, A., Fuzzy group theory, Springer Verlag, Berlin, 2005.
\bibitem{7} V. Murali, B.B. Makamba, On an equivalence of fuzzy subgroups, I, Fuzzy Sets and Systems {\bf 123} (2001), 259-264.
\bibitem{8} V. Murali, B.B. Makamba, On an equivalence of fuzzy subgroups, II, Fuzzy Sets and Systems {\bf 136} (2003), 93-104.
\bibitem{9} V. Murali, B.B. Makamba, On an equivalence of fuzzy subgroups, III, Int. J. Math. Sci. {\bf 36} (2003), 2303-2313.
\bibitem{10} V. Murali, B.B. Makamba, Counting the number of fuzzy subgroups of an abelian group of order $p^n q^m$, Fuzzy Sets and Systems {\bf 144} (2004), 459-470.
\bibitem{11} V. Murali, B.B. Makamba, Fuzzy subgroups of finite abelian groups, FJMS {\bf 14} (2004), 113-125.
\bibitem{12} R. Schmidt, Subgroup lattices of groups, de Gruyter Expositions in Mathematics 14, de Gruyter, Berlin, 1994.
\bibitem{13} R.P. Stanley, Enumerative combinatorics, II, Cambridge University Press, Cambridge, 1999.
\bibitem{14} M. Suzuki, Group theory, I, II, Springer Verlag, Berlin, 1982, 1986.
\bibitem{15} M. T\u arn\u auceanu, Groups determined by posets of subgroups, Ed. Matrix Rom, Bucure\c sti, 2006.
\bibitem{16} M. T\u arn\u auceanu, The number of fuzzy subgroups of finite cyclic groups and Delannoy numbers, European J. Combin. {\bf 30} (2009), 283-287, doi: 10.1016/j.ejc.2007.12.005.
\bibitem{17} M. T\u arn\u auceanu, Distributivity in lattices of fuzzy subgroups, Inform. Sci. {\bf 179} (2009), 1163-1168, doi: 10.1016/j.ins.2008.12.003.
\bibitem{18} M. T\u arn\u auceanu, Classifying fuzzy subgroups of finite nonabelian groups, Iran. J. Fuzzy Syst. {\bf 9} (2012), 33-43.
\bibitem{19} M. T\u arn\u auceanu, A note on the lattice of fuzzy subgroups of a finite group, J. Mult.-Valued Logic Soft Comput. {\bf 19} (2012), 537-545.
\bibitem{20} M. T\u arn\u auceanu, On the number of fuzzy subgroups of finite symmetric groups, J. Mult.-Valued Logic Soft Comput. {\bf 21} (2013), 201-213.
\bibitem{21} M. T\u arn\u auceanu, Classifying fuzzy subgroups for a class of finite $p$-groups, Critical Review {\bf 7} (2013), 30-39.
\bibitem{22} M. T\u arn\u auceanu, Classifying fuzzy normal subgroups of finite groups, Iran. J. Fuzzy Syst. {\bf 12} (2015), 107-115.
\bibitem{23} M. T\u arn\u auceanu, The number of chains of subgroups of a finite elementary abelian $p$-group, to appear in  Sci. Bull., Ser. A, Appl. Math. Phys., Politeh. Univ. Buchar., 2015.
\bibitem{24} M. T\u arn\u auceanu, L. Bentea, On the number of fuzzy subgroups of finite abelian groups, Fuzzy Sets and Systems {\bf 159} (2008), 1084-1096, doi: 10.1016/j.fss.2007.11.014.
\bibitem{25} M. T\u arn\u auceanu, L. Bentea, A note on the number of fuzzy subgroups of finite groups, Sci. An. Univ. "Al.I. Cuza" Ia\c si, Math., {\bf 54} (2008), 209-220.
\bibitem{26} I. Tomescu, Introduction to combinatorics, Collet's Publishers Ltd., London, 1975.
\bibitem{27} A.C. Volf, Counting fuzzy subgroups and chains of subgroups, Fuzzy Systems \& Artificial Intelligence {\bf 10} (2004), 191-200.
\bibitem{28} A. Weinberger, Reducing fuzzy algebra to classical algebra, New Math. Natur. Comput. {\bf 1} (2005), 27-64.
\end{thebibliography}
\end{document}